\documentstyle{article}

\newcommand{\be}{\begin{equation}}
\newcommand{\ee}{\end{equation}}
\newcommand{\ba}{\begin{eqnarray}}
\newcommand{\ea}{\end{eqnarray}}
\newcommand{\baa}{\begin{eqnarray*}}
\newcommand{\eaa}{\end{eqnarray*}}
\newcommand{\bb}{\begin{thebibliography}}
\newcommand{\eb}{\end{thebibliography}}

\newcommand{\bi}[1]{\bibitem{#1}}
\newcommand{\lab}[1]{\label{#1}}
\newcommand{\re}[1]{(\ref{#1})}

\newcounter{my}
\newcommand{\he}%
   {\stepcounter{equation}\setcounter{my}%
   {\value{equation}}\setcounter{equation}0%
   }%
\newcommand{\she}%
   {\setcounter{equation}{\value{my}}%
    }%

\renewcommand\t{\tilde}

\newcommand\ve{\varepsilon}

\newtheorem{pr}{Proposition}

\begin{document}

\vspace*{10mm}

\begin{center}

{\Large \bf A ``missing`` family of classical orthogonal polynomials}

\vspace{5mm}

{\large \bf Luc Vinet}

\medskip

{\it Centre de recherches math\'ematiques
Universite de Montr\'eal, P.O. Box 6128,
Centre-ville Station, Montr\'eal (Qu\'ebec), H3C 3J7}

\medskip

and

\medskip

{\large \bf Alexei Zhedanov}

\medskip

{\em Donetsk Institute for Physics and Technology, Donetsk 83114,  Ukraine}

\end{center}


\begin{abstract}
We study a family of "classical" orthogonal polynomials which satisfy (apart from a 3-term recurrence relation)
an eigenvalue problem with a differential operator of Dunkl-type.
These polynomials can be obtained from the little $q$-Jacobi polynomials in the limit $q=-1$.
We also show that these polynomials provide a nontrivial realization of the Askey-Wilson algebra for $q=-1$.

\vspace{2cm}


\vspace{2cm}


\end{abstract}

\bigskip\bigskip


\newpage
\section{Introduction}
The Askey scheme \cite{KS}, \cite{KLS} provides a list of all known "classical'' orthogonal polynomials. The term ``classical" means that apart from  a three-term recurrence relation
\be
P_{n+1}(x) + b_n P_n(x) + u_n P_{n-1}(x) = x P_n(x), \quad P_{-1}=0, \; P_0=1 \lab{3term_P} \ee
(which is a general feature of orthogonal polynomials \cite{Chi}, \cite{Sz}) these polynomials satisfy also an eigenvalue equation
\be L P_n(x) = \lambda_n P_n(x) \lab{LPn} \ee
where the operator $L$ acts on the variable $x$. The operator $L$ can be either a second-order differential operator (in this case we obtain the purely classical polynomials: Jacobi, Laguerre, Hermite) or a second-order difference operator on a uniform or non-uniform grid. The latter case leads to Wilson ($q$=1) and Askey-Wilson ($q \ne 1$) polynomials and their descendents.

Bannai and Ito \cite{BI} proved a general theorem (which generalizes a famous theorem of Leonard \cite{Leonard}, dealing only with polynomials orthogonal on a finite set of points) stating that all such orthogonal polynomials should coincide with the Askey-Wilson polynomials or their specializations. They also found a ''missing" case of the Askey scheme corresponding to the limit $q=-1$ of the $q$-Racah polynomials (see also \cite{Ter}). The example of Bannai and Ito corresponds to polynomials orthogonal on a finite set of points. Hence in this case the operator $L$ is merely a finite-dimensional 3-diagonal matrix which corresponds to the case of "Leonard pairs" (see, e.g. \cite{Ter}).

It is hence sensible to investigate other possibilities as $q$ approaches -1 in the Askey scheme. Of course, the limit $q=1$ is well studied and classified (see, e.g. \cite{KS}). The limit $q=-1$ however, has not been explored much. In \cite{AI} Askey and Ismail  have studied the limit $q=-1$ for the $q$-ultraspherical polynomials, but in this case, the operator $L$ disappears in the limit $q=-1$.

Here we show that there is a very simple class of polynomials which can be obtained from the little $q$-Jacobi polynomials in the limit $q=-1$. Under appropriate choice of the parameters, the operator $L$ survives in the limit $q=-1$. The polynomials thus obtained are indeed classical: they satisfy the eigenvalue equation \re{LPn}. But in contrast to the case of pure classical polynomials (like Jacobi polynomials), the operator $L$ is a combination of a differential operator of first order and of the reflection operator $R$ (see formula \re{L_0_R} in the next section). Operators of this type are known as Dunkl operators \cite{Dunkl}. So far operators of Dunkl type were used typically to transform one polynomial family into another. In the one-dimensional case, the Dunkl operator $T_{\mu}$ is defined as  \cite{Dunkl}, \cite{Cheikh}
\be
T_{\mu}f(x) =  f'(x) + \mu \frac{f(x)-f(-x)}{x} \lab{Dunkl_T} \ee
where $\mu$ is a deformation parameter (when $\mu=0$ the Dunkl operator becomes the ordinary derivative operator).

Clearly, the Dunkl operator \re{Dunkl_T} reduces the degree of any polynomial by 1. Hence, there are no polynomials which are eigenfunctions of the operator $T_{\mu}$.

Our operator \re{LPn}, in contrast, preserves the linear space of polynomials of given degree. In the following we construct the general solution of the eigenvalue equation \re{L_0_P} and obtain explicit expression for the corresponding little -1 Jacobi polynomials $P_n^{(-1)}(x)$.

We also show that the polynomials $P_n^{(-1)}(x)$ are Dunkl-classical, i.e. $T_{\mu}P_n^{(-1)}(x)=Q_{n-1}(x)$, where $Q_n(x)$ is another set of little -1 Jacobi polynomials (with different parameters)

\section{Limit of the little $q$-Jacobi polynomials as $q \to -1$}
\setcounter{equation}{0} The little $q$-Jacobi polynomials are
defined through the recurrence coefficients \be u_n = A_{n-1}C_n,
\quad b_n = A_n + C_n, \lab{ub_lqJ} \ee where $A_n, C_n$ are given
by
$$
A_n = q^n
\frac{(1-aq^{n+1})(1-abq^{n+1})}{(1-abq^{2n+1})(1-abq^{2n+2})}, \;
C_n = aq^n \frac{(1-q^{n})(1-bq^{n})}{(1-abq^{2n+1})(1-abq^{2n})}.
$$
They have the following simple expression in terms of the basic hypergeometric
function \be P_n(x) = \kappa_n \; {_2}\Phi_1\left({q^{-n},
abq^{n+1} \atop aq} \Bigl | q; qx  \right) \lab{lqJ_hyp} \ee with
a normalization factor $\kappa_n$ to ensure that they are monic.

They satisfy the orthogonality relation \be \sum_{s=0}^{\infty}
\frac{(bq;q)_s}{(q;q)_s} (aq)^s P_n(q^s)P_m(q^s) = h_n \:
\delta_{nm}. \lab{lqJ_ort} \ee

The moments corresponding to this weight function are \be
c_n=\frac{(aq;q)_n}{(abq^2;q)_n}, \lab{lqj_moms} \ee where
$(x;q)_n = (1-x)(1-qx) \dots(1-q^{n-1}x)$ is standard notation for
the $q$-shifted factorials (Pochhammer $q$-symbol).

There is a $q$-difference equation of the form \be
a(bq-x^{-1})(P_n(qx)-P_n(x)) +(1-x^{-1})(P_n(q^{-1}x)-P_n(x)) =
\lambda_n P_n(x), \lab{lqJ_eq} \ee where \be \lambda_n =
(q^{-n}-1)(1-abq^{n+1}). \lab{lqJ_lam} \ee If $a=q^{\alpha},
b=q^{\beta}$ then in the limit $q\to 1$ we get the ordinary Jacobi
polynomials with parameters $\alpha, \beta$.

There is, however, another nontrivial limit if one puts \be
q=-e^{\ve}, \; a=-e^{\ve \alpha}, \; b=-e^{\ve \beta} \lab{qab_e}
\ee and take the limit $\ve \to 0$. This is the limit $q=-1$ of the little $q$-Jacobi polynomials.

A direct calculation shows that in this limit, we have the recurrence
coefficients \be u_{n} = \frac{(n+(1-\theta_n)\alpha)(n+\beta+
\theta_n \alpha)}{(2n+\alpha+\beta)^2}, \quad b_n = (-1)^n \:
\frac{(2n+1)\alpha + \alpha \beta +\alpha^2 +(-1)^n
\beta}{(2n+\alpha+\beta)(2n+2+\alpha+\beta)}, \lab{lqJ-1_ub} \ee
where
$$
\theta_n= \frac{1+(-1)^n}{2}
$$
is the characteristic function of even numbers.

The corresponding moments are obtained directly from the moments \re{lqj_moms}: \be
c_{2n}=c_{2n-1}=\frac{(\alpha/2+1/2)_n}{(\alpha/2+\beta/2+1)_n},
\quad n=1,2,3,\dots,  \lab{c_lim} \ee
where $(x)_n=x(x+1) \dots(x+n-1)$ is the ordinary Pochhammer symbol (shifted factorial).

Using this explicit expression for the moments, we can recover the weight function $w(x)$ for the resulting orthogonal polynomials.
It is easily verified that
\be w(x)=\kappa |x|^{\alpha}
(1-x^2)^{(\beta-1)/2}(1+x), \lab{weight} \ee where
$$
\kappa= \frac{\Gamma(\alpha/2+\beta/2+1)}{\Gamma(\beta/2+1/2)
\Gamma(\alpha/2+1/2)}.
$$
Indeed, we have (using the ordinary Euler $B$-integral)
$$
\int_{-1}^1 w(x) x^n dx = c_n, \quad n=0,1,2,\dots,
$$
where $c_n$ is given by \re{c_lim}. The coefficient $\kappa$ is
chosen to provide the standard normalization condition $c_0=1$.
Under the obvious conditions $\alpha>-1, \: \beta>-1$, the weight
function $w(x)$ is positive, all moments $c_n$ are well defined
and the moment problem is positive definite, i.e. $\Delta_n >0$
for all $n=0,1,2,\dots$.

Consider the form of the $q$-difference equation
\re{lqJ_eq} in this limit. We divide both sides of \re{lqJ_eq} by $\ve$ and
introduce the operator $L_{\ve}$ which acts on any polynomial
$f(x)$ as \be L_{\ve} f(x) = a\ve^{-1}(bq-x^{-1})(f(qx)-f(x))
+\ve^{-1}(1-x^{-1})(f(q^{-1}x)-f(x)) \lab{lqJ_L} \ee (the
parameters $q,a,b$ depends on $\ve$ as in \re{qab_e}). For
monomials $f=x^n$ we have in the limit $\ve =0$ \be L_0 x^n =
\xi_n x^n + \eta_n x^{n-1}, \lab{Lxn} \ee where \be \xi_n =
2(-1)^{n+1} n +(1-(-1)^n)(\alpha+\beta+1), \quad \eta_n =2(-1)^n n
-(1-(-1)^n)\alpha. \lab{xi_eta_lim} \ee This allows one to present
the operator $L_0$ in the form \be L_0=2(1-x)\partial_x R
+(\alpha+\beta+1-\alpha x^{-1})(1-R), \lab{L_0_R} \ee where $R$ is
the reflection operator $Rf(x)=f(-x)$.

Thus we have that our polynomials are classical: they satisfy the eigenvalue equation
\be
L_0 P_n(x) = \lambda_n P_n(x), \lab{L_0_P} \ee
where
\be
\lambda_n = \left\{ -2n \quad \mbox{if $n$}  \; \mbox{is even} \atop 2(\alpha+\beta+n+1) \quad \mbox{if $n$}  \; \mbox{is odd} .  \right .
\lab{lambda_n_eo} \ee
But in contrast to true classical polynomial the operator $L_0$ is not purely differential: it contains the reflection operator $R$.

In the next section we construct  the general solution of the eigenvalue problem \re{L_0_P} in terms of the Gauss hypergeometric functions.

\section{Different forms of $R$-differential equation}
\setcounter{equation}{0} Consider the eigenvalue equation \be L_0
F(x) = \lambda F(x), \lab{eiG_L_0} \ee where $L_0$ has the
expression \re{L_0_R}. The operator $L_0$ is a differential
operator of the first order containing the reflection operator
$R$. Notice the important property of the operator $L_0$: it preserves any linear space of polynomials of degrees $\le N$ for any $N=1,2,3,\dots$.
Hence the operator $L_0$ behaves like the classical hypergeometric operator: for any $n=0,1,2, \dots$ there exists a polynomial eigenvalue solution
$L_0 P_n(x) = \lambda_n P_n (x)$ for an appropriate set of eigenvalues $\lambda_n$.

This solution can be found directly from the property \re{Lxn} of the operator $L_0$. Indeed, in the basis $x^n$ the operator $L_0$ is two-diagonal and hence the corresponding eigenvalue solutions can be found explicitly because the coefficients $\xi_n$ and $\eta_n$ are known.

We consider here another method, which, though more complicated, demonstrates the relation with the Gauss hypergeometric equation.

We can arrive at a pure differential equation (without the
operator $R$) by a standard procedure. Let us present the
function $F(x)$ as a superposition of the even and odd parts:
$$
F(x) = f(x) + g(x),
$$
where $f(-x)=f(x), \; g(-x)=-g(x)$. The functions $f(x)$ and
$g(x)$ are determined uniquely from $F(x)$. The operator $R$ acts
on these functions as: $R f(x)=f(x), \; R g(x) =- g(x)$. This
allows one to rewrite the eigenvalue equation \re{eiG_L_0} in the form
\be (1-x) f'(x) + (x-1) g'(x) + (1+\alpha+\beta - \alpha/x)g(x) =
\lambda(f(x)+g(x))/2. \lab{eigen_fG_1} \ee Consider also the associated
eigenvalue equation \be R L_0 F(x) = \lambda R F(x) \lab{eiG_RL_0}
\ee which is obtained from \re{eiG_L_0} by application of the
operator $R$. This equation can be presented in the form \be (1+x)
f'(x) + (1+x) g'(x) + (1+\alpha+\beta + \alpha/x)g(x) =
\lambda(g(x)-f(x))/2 \lab{eigen_fG_2} \ee Adding and subtracting,
we obtain the simpler system of equations \ba &&f'(x) + xg'(x) +
(1+\alpha+\beta)g(x) =\lambda g(x) /2, \nonumber \\&&xf'(x) +
g'(x) +\alpha g(x)/x = -\lambda f(x)/2 \lab{sys_fg} \ea This
casts the eigenvalue equation \re{eigen_fG_1} equivalently
as a system of two linear differential equations of
first order. We can eliminate the function $g(x)$ from this
system:
\be
g(x) = \frac{2(x^2-1) f'(x) + \lambda x f(x)}{2(\beta+1) -\lambda}
\lab{g_f} \ee
and obtain a second-order differential equation for $f(x)$: \be 4
x(x^2-1) f''(x) +  4((\alpha+\beta+3)x^2-\alpha) f'(x)  + \lambda
x(2(\alpha+\beta) + 4 -\lambda)f(x) =0. \lab{eq_f} \ee By a change
of variable $x \to x^2$ one can reduce equation \re{eq_f} to
the Gauss hypergeometric equation thus obtaining the general
solution: \be f(x) = C_1 \: {_2}F_1 \left( {\frac{\lambda}{4},
\frac{\alpha+\beta}{2}+1 -\frac{\lambda}{4} \atop
\frac{\alpha+1}{2} }; x^2  \right) + C_2 \: x^{1-\alpha}{_2}F_1
\left( {\frac{\lambda+2-2\alpha}{4}, \frac{2\beta+6-\lambda}{4}
\atop \frac{3-\alpha}{2} }; x^2  \right), \ee where $C_1, C_2$ are
arbitrary constants. By construction, the function $f(x)$ must be even
for all values of the parameters $\alpha, \beta$. This is possible
only if $C_2 =0$ and the solution for $f(x)$ contains
only one undetermined constant: \be f(x) = C_1 \: {_2}F_1 \left(
{\frac{\lambda}{4}, \frac{\alpha+\beta}{2}+1 -\frac{\lambda}{4}
\atop \frac{\alpha+1}{2} }; x^2  \right). \lab{sol_f} \ee

From formula \re{g_f} we obtain the solution for the function $g(x)$
\be g(x) = -\frac{\lambda C_1}{2(\alpha+1)} \: x \:  {_2}F_1 \left( {1+\frac{\lambda}{4},
\frac{\alpha+\beta}{2}+1 -\frac{\lambda}{4} \atop
\frac{\alpha+3}{2} }; x^2  \right). \lab{sol_g} \ee
It is clear that the function $g(x)$ is odd as desired.

In formula \re{g_f} it was assumed that $\lambda \ne 2(\beta+1)$. If $\lambda = 2(\beta+1)$, then the Gauss hypergeometric function ${_2}F_1(x^2)$ reduces to the elementary function ${_1}F_0(x^2)$ and for the function $f(x)$ we obtain
\be
f(x) = C_1 \: (1-x^2)^{-\frac{\beta+1}{2}}.
\lab{f_elem} \ee
Similarly, if $\lambda=2(\beta-1)$ we obtain an elementary expression
$$
g(x) = -\frac{C_1(\beta-1)}{\alpha+1} \: (1-x^2)^{-\frac{\beta+1}{2}}
$$
for the function $g(x)$.

Now we have the
general solution of equation \re{eiG_L_0} in the form \be F(x) =
f(x)+g(x) = C(\lambda) \: {_2}F_1 \left( {\frac{\lambda}{4},
\frac{\alpha+\beta}{2}+1 -\frac{\lambda}{4} \atop
\frac{\alpha+1}{2} }; x^2  \right) - \frac{\lambda C(\lambda)}{2(\alpha+1)} \: x \: {_2}F_1
\left( {1+\frac{\lambda}{4}, \frac{\alpha+\beta}{2}+1
-\frac{\lambda}{4} \atop \frac{\alpha+3}{2} }; x^2  \right),
\lab{gen_sol_F1} \ee where the coefficient $C(\lambda)$
may depend on $\lambda$ (as well as the parameters $\alpha,
\beta$).

We would like to get polynomial solutions, i.e. find
eigenvalues $\lambda_n, \; n=0,1,2,\dots$ such that $F_n(x)$ is a
polynomial in $x$ of exact degree $n$. It is easily seen that
polynomial solutions are possible only if either $\lambda=-4n, \:
n=0,1,2,\dots$ or $\lambda=2(\alpha+\beta+2+2n), \;
n=0,1,2,\dots$. If $\lambda=-4n$, the first term in
\re{gen_sol_F1} is a polynomial of degree $2n$ whereas the second
term is a polynomial of degree $2n-1$. Hence for all $\lambda_n =
-4n$ we will have polynomials $F_n(x)$ of the even degree $2n$. If
$\lambda=2(\alpha+\beta+2+2n)$ then the first term in
\re{gen_sol_F1} is a polynomial of degree $2n$, while the second
term is a polynomial of degree $2n+1$. Hence for all
$\lambda_n=2(\alpha+\beta+2+2n), \; n=0,1,2,\dots$ we obtain
polynomials $F_{n}(x)$ having odd degree $2n+1$. This solves the problem,
and we thus have the following explicit expressions.

If $n$ is even then
\be
P_n^{(-1)}(x)= \kappa_n \left[\: {_2}F_1 \left({-\frac{n}{2}, \frac{n+\alpha+\beta+2}{2} \atop  \frac{\alpha+1}{2}}; x^2  \right)  +  \frac{nx}{\alpha+1} \: {_2}F_1 \left({1-\frac{n}{2}, \frac{n+\alpha+\beta+2}{2} \atop  \frac{\alpha+3}{2}}; x^2  \right) \right] \lab{eve_P} \ee

If $n$ is odd then
\be
P_n^{(-1)}(x)= \kappa_n \left[\: {_2}F_1 \left({\frac{1-n}{2}, \frac{n+\alpha+\beta+1}{2} \atop  \frac{\alpha+1}{2}}; x^2  \right) \lab{even_P} - \frac{(\alpha+\beta+n+1)x}{\alpha+1} \: {_2}F_1 \left({\frac{1-n}{2}, \frac{n+\alpha+\beta+3}{2} \atop  \frac{\alpha+3}{2}}; x^2  \right) \right], \lab{odd_P} \ee
where $\kappa_n$ is an appropriate normalization factor to ensure that the polynomial $P_n^{(1)}(x)=x^n + O(x^{n-1})$ is monic.
(We need not its explicit expression).

\section{Relation with the symmetric Jacobi polynomials}
\setcounter{equation}{0}
The symmetric Jacobi polynomials were introduced by Chihara \cite{Chi}. They can be defined as follows.

Let $P_n^{(\xi,\eta)}(x)$ be the (modified) Jacobi polynomials \be
P_n^{(\xi,\eta)}(x) = \kappa_n \: {_2}F_1 \left({-n, n+\xi+\eta +1
\atop \xi+1};x  \right) \lab{P_J} \ee which are orthogonal on the
interval $[0,1]$
$$
\int_0^1 P_n^{(\xi,\eta)}(x) P_m^{(\xi,\eta)}(x) x^{\xi}(1-x)^\eta
dx = h_n \delta_{nm}
$$
($\kappa_n$ is a factor needed for the polynomials $P_n^{\xi,\eta}(x)$ to be monic).

We can introduce symmetric polynomials $S_n^{(\xi,\eta)}(x)$ by
the formulas \be S_{2n}^{(\xi,\eta)}(x) = P_n^{(\xi,\eta)}(x^2),
\quad S_{2n+1}^{(\xi,\eta)}(x) = xP_n^{(\xi+1,\eta)}(x^2).
\lab{S_P} \ee Then it is easily verified that the polynomials
$S_n(x)$ satisfy the symmetric property $S_n(-x)=(-1)^n S_n(x)$
and are orthogonal on the interval $[-1,1]$ \be \int_{-1}^1
S_n^{(\xi,\eta)}(x) S_m^{(\xi,\eta)}(x)w(x) dx= h_n \delta_{nm}
\lab{ort_S} \ee with the weight function \be w(x)=
|x|^{2\xi+1}(1-x^2)^{\eta}. \lab{w_sym_J} \ee (Note that the
original definition  \cite{Chi} by Chihara of the symmetric Jacobi
polynomials $S_n(x)$ is slightly different. This is not essential
for our purposes).

If $\xi=-1/2$, the polynomials $S_m^{(-1/2,\eta)}(x)$ have the weight function $w(x)=(1-x^2)^{\eta}$ and hence can be identified with the classical ultraspherical (Gegenbauer) polynomials. This is why the polynomials $S_n^{\xi,\eta}(x)$ are  sometimes called generalized Gegenbauer polynomials \cite{Belm}.

Starting from the polynomials $S_n^{\xi,\eta}(x)$, we can perform the Christoffel transform \cite{SVZ}
\be
\t S_n^{\xi,\eta}(x) = \frac{S_{n+1}^{\xi,\eta}(x)  - A_n S_n^{\xi,\eta}(x)}{x+1}, \lab{CT_S} \ee
where
$$
A_n = \frac{S_{n+1}^{\xi,\eta}(-1)}{S_n^{\xi,\eta}(-1)}.
$$
The polynomials $\t S_n^{\xi,\eta}(x)$ are again orthogonal on the interval $[-1,1]$ with the weight function $\t w(x)$ obtained from \re{w_sym_J} as
\be
\t w(x) = w(x) (x+1) = |x|^{2\xi+1}(1-x^2)^{\eta}(1+x) \lab{w_t} \ee
(we do not take into account the normalization condition which is not essential for our purposes).

The Christoffel transform \re{CT_S} can be presented in explicit form because the coefficients $A_n$ are expressible in terms of the hypergeometric function ${_2}F_1 (z)$ with argument $z=1$ (this is seen from formulas \re{S_P} and \re{P_J}).

Hence from \re{CT_S}, we have an explicit expression for the
polynomials $\t S_n(x)$ as a linear combination of two
hypergeometric functions. The polynomials $\t S_n(x)$ are not
symmetric, i.e. they satisfy the general 3-term recurrence
relation \be \t S_{n+1}(x) + \t b_n \t S_{n}(x) + \t u_n \t
S_{n-1}(x) = x \t S_{n}(x), \lab{3-term_t_S} \ee where the
coefficients $\t b_n, \: \t u_n$ can be expressed in terms of the
coefficients $A_n$ using the properties of the Christoffel
transform  \cite{SVZ}. Hence these coefficients are also explicit.

One can express the weight function \re{w_t} in an equivalent form
$$
\t w(x) = |x|^{2\xi+1}(1-x^2)^{\eta+1}(1-x)^{-1}
$$
which corresponds to the Geronimus transformation \cite{Ger1}, \cite{ZhS}:
\be
\t S_n(x)  =  S_n^{\xi,\eta+1}(x) - B_n S_{n-1}^{\xi,\eta+1}(x) \lab{S_Ger} \ee
with some coefficients $B_n$.

We thus have that the same polynomials $\t S_n(x)$ can be obtained  from the generalized Gegenbauer polynomials $S_n^{\xi,\eta}(x)$ by either Christoffel or Geronimus transformations.

In \cite{Chi_Chi} these polynomials were presented in the form \re{S_Ger}; the recurrence coefficients $\t b_n, \: \t u_n$ were derived as well in \cite{Chi_Chi}.

The comparison of the weight functions \re{w_t} and \re{weight} leads to the conclusion that the little -1 Jacobi polynomials   $P_n^{(-1)}(x)$ coincide with polynomials $\t S_n^{(\xi, \eta)}$, where
$$
\xi=\frac{\alpha-1}{2}, \quad \eta=\frac{\beta-1}{2}.
$$
We thus identified our "classical" polynomials $P_n^{(-1)}(x)$ with the "nonsymmetric" generalized Gegenbauer proposed by L.Chihara and T.Chihara \cite{Chi_Chi}.

Explicitly we have from \re{S_Ger}
\be
P_n^{(-1)}(x) = S_n^{(\frac{\alpha-1}{2},\frac{\beta-1}{2})} (x) - B_n S_{n-1}^{(\frac{\alpha-1}{2},\frac{\beta-1}{2})}(x), \lab{litll_P_S} \ee
where
$$
B_n= \frac{2n+(1-(-1)^n) \alpha}{2(\alpha+\beta+2n)}.
$$

The same polynomials were also considered in \cite{Atia}, \cite{Atia2} from another point of view.

We would like to stress that other basic properties of these polynomials , e.g. the existence of a Dunkl-type operator $L_0$ providing the eigenvalue problem \re{L_0_P}, as well as the origin of these polynomials as a limit case ($q$=-1) of the little $q$-Jacobi polynomials had not been identified.

\section{The Dunkl-classical property}
\setcounter{equation}{0}
All "classical" orthogonal polynomials satisfy an important characteristic condition: they are "covariant" with respect to a "derivative" operator $\cal D$:
\be
{\cal D} P_n(x) = [n] \t P_{n-1}(x), \lab{DPP} \ee
where $[n]$ is a specific function of $n$ depending on the choice of the operator $\cal D$, $\t P_n(x)$ is another set of "classical" orthogonal polynomials, and the operator $D$ possesses the basic property of reducing the degree any polynomial by one.

Well known examples of the operator ${\cal D}$ are:

(i) the derivative operator ${\cal D} = \partial_x$;

(ii) the difference operator ${\cal D} f(x) = f(x+1)-f(x)$;

(iii) the $q$-derivative operator: ${\cal D} f(x) = \frac{f(xq)-f(x)}{(q-1)x}$

(iv) the Askey-Wilson operator ${\cal D} f(x(s)) = \frac{f(x(s+1/2)) - f(x(s-1/2))}{x(s+1/2) -x(s-1/2)}$.

In the last case the function $x(s)$ is either trigonometric $x(s) = a_1 q^s + a_2 q^{-s} + a_0$ or quadratic $x(s)=a_2 s^2 + a_1 s + a_0$ with some constants $a_i$.

Recently it was recognized that apart from these operators there is one more operator which generates "classical" orthogonal polynomials. This operator is the Dunkl operator $T_{\mu}$ \re{Dunkl_T}. Namely, in \cite{Cheikh} it was shown that the only {\it symmetric} orthogonal polynomials $P_n(x)$ satisfying the property
\be
T_{\mu} P_n(x) = [n]_{\mu} \t P_{n-1}(x), \quad [n]_{\mu} = n + (1-(-1)^n)\mu \lab{TPP} \ee
are the generalized Hermite or the generalized Gegenbauer polynomials. Recall that symmetric orthogonal polynomials are defined by the property $P_n(-x) = (-1)^n P_n(x)$. The generalized Hermite polynomials $H_n^{(\mu)}(x)$ \cite{Chi} are symmetric orthogonal polynomials which are orthogonal on whole real line with the weight function
$$
w(x)= |x|^{2 \mu} \exp(-x^2)
$$
When $\mu=0$ (i.e. in the case when the Dunkl operator $T_{\mu}$ becomes the derivative operator $\partial_x$) the generalized Hermite polynomials become the ordinary Hermite polynomials.

The generalized Gegenbauer polynomials $S_n^{(\xi, \eta)}(x)$ \cite{Chi}, \cite{Belm} are orthogonal on the interval $[-1,1]$ with the weight function \re{w_sym_J}. The generalized Gegenbauer polynomials satisfy the Dunkl-classical property \re{TPP} with $\mu=\xi+1/2$

In both cases it is assumed that $\mu>-1/2$.

In the present case we have correspondingly  the following simple but important result
\begin{pr}
The little -1 Jacobi polynomials $P_n^{(-1)}(x)$  satisfy the Dunkl-classical property \re{TPP} with $\mu=\alpha/2$, where the polynomials $\t P_n(x)$  are again little -1 Jacobi polynomials with parameters $(\alpha, \beta+2)$.
\end{pr}
The proof of this proposition follows easily from the explicit formula \re{litll_P_S} and from the fact that the generalized Gegenbauer polynomials $S_n^{(\xi,\eta)}(x)$ satisfy the Dunkl-classical propertry  \re{TPP} \cite{Cheikh}:
$$
T_{\mu} S_n^{(\xi,\eta)}(x) = [n]_{\mu} S_{n-1}^{(\xi,\eta+1)}(x), \quad \mu=\xi+1/2
$$
In contrast to the assumptions of \cite{Cheikh}, the little -1 Jacobi polynomials are not symmetric.
Hence we perhaps obtained  the first example of Dunkl-classical orthogonal polynomials beyond the family of symmetric polynomials.
The problem of finding all such orthogonal polynomials is an interesting open question.

All known families of "classical" orthogonal polynomials possess
not only lowering operators like \re{DPP} but also raising
operators $\Theta$ with the property \be \Theta P_n(x) = \nu_{n+1}
Q_{n+1}(x), \lab{raising_P} \ee where the polynomials $Q_n(x)$
belong to the same family of classical orthogonal polynomials
(albeit with different parameters).

 In the case of the little -1 Jacobi polynomials it is directly verified that the operator $\Theta$ does exists and has the expression
 \be
 \Theta f(x) = (x^2-1) f'(x) + \frac{\alpha (x-1)^2}{2x} f(-x) + \left((\beta+\alpha/2)x -1-\frac{\alpha}{2x} \right) f(x) \lab{raising_R} \ee
Given \re{raising_R}, property \re{raising_P} holds with
$$
\nu_{n} = n+\beta-1 +  \frac{1-(-1)^n}{2} \alpha = \beta-1+[n]_{\mu}, \quad \mu=\alpha/2
$$
and $Q_n(x)$ the same monic little -1 Jacobi polynomials with parameters $(\alpha,\beta-2)$ .

The generalized Hermite and Gegenbauer polynomials can be obtained from the ordinary Hermite and Gegenbauer polynomials through the acting of the
Dunkl intertwining operator \cite{Dunkl}, \cite{Cheikh}.

Recall that the Dunkl intertwining operator $V_{\mu}$ acts on the
space of polynomials by the formulas \cite{Dunkl}, \cite{Dunkl2}
$$
V_{\mu} x^n = \sigma_n x^n, \quad \sigma_{2n-1}=\sigma_{2n} = \frac{(1/2)_n}{(\mu+1/2)_n}, $$
and is realized by the following integral representation \cite{Dunkl}
$$
V_{\mu} (f(x)) = \frac{\Gamma(\mu+1/2)}{\Gamma(\mu) \Gamma(1/2)}
\: \int_{-1}^1 f(xt) (1-t)^{\mu-1}(1+t)^{\mu} dt
$$

It preserves the space of polynomials and has the fundamental
intertwining property \be T_{\mu} V_{\mu}= V_{\mu} \partial_x.
\lab{inter_V} \ee From this property it is possible to obtain the
following result (see Proposition {\bf 2} below). Assume that the
monic polynomials $P_n(x)$ and $Q_n(x)$ are related as \be P_n'(x)
= n Q_{n-1}(x). \lab{DPQ} \ee Let us construct the monic
polynomials $\t P_n(x) = \sigma_n^{-1}\: V_{\mu} P_n(x), \; \t
Q_n(x) = \sigma_n^{-1}\: V_{\mu} Q_n(x)$. Then these polynomials
are correspondingly related: \be T_{\mu} \t P_n(x) = [n]_{\mu} \t
Q_{n-1}(x). \lab{TPTQ} \ee In particular, if all polynomials
$P_n(x), Q_n(x), \t P_n(x), \t Q_n(x)$ are orthogonal then the
operator $V_{\mu}$ allows to obtain Dunkl-classical polynomials
(defined by property \re{TPTQ}) from ordinary classical
polynomials (defined by property \re{DPQ}).

In \cite{Cheikh} it was shown that the generalized Gegenbauer
polynomials  $S_{n}^{(\xi,\eta)}(x)$ can be obtained from the
ordinary Gegenbauer polynomials $S_{n}^{(-1/2,\eta+\mu)}(x)$ by
the action of the intertwining operator $V_{\mu}$: \be
S_{n}^{(\xi,\eta)}(x) = \sigma_n^{-1}V_{\mu}
S_{n}^{(-1/2,\eta+\mu)}(x), \quad \mu=\xi+1/2. \lab{Ch_S} \ee (A
similar property for the generalized Hermite polynomials was
obtained by Dunkl \cite{Dunkl}, \cite{Dunkl2}).

Introduce now the ordinary monic Jacobi polynomials $P_n^{(\xi,\eta)}(x)$ by the formula
\be
P_n^{(\xi,\eta)}(x) = \frac{2^n (\xi+1)_n}{(\xi+\eta+n+1)_n} \: {_2}F_1 \left({-n, n+\xi+\eta +1 \atop \xi+1}; \frac{1-x}{2}  \right).
\lab{KS_PJ} \ee
Notice that this definition differs from \re{P_J} by an affine transformation of the argument. The polynomials \re{KS_PJ} coincide with standard Jacobi polynomials orthogonal on the interval $[-1,1] $\cite{KS}.

We have the following
\begin{pr}
The little -1 Jacobi polynomials \re{eve_P}, \re{odd_P} can be obtained from the Jacobi polynomials \re{KS_PJ} by the action of the
Dunkl intertwining operator
\be
P_n^{(-1)}(x) = \sigma_n^{-1} V_{\mu} P_n^{(\xi,\xi+1)}(x), \quad \xi=\frac{\alpha+\beta-1}{2}, \; \mu = \frac{\alpha}{2}. \lab{PVP} \ee
\end{pr}
The proof of this proposition is based on formula \re{litll_P_S} and property \re{Ch_S}.

\section{Askey-Wilson algebra relations for the little -1 Jacobi polynomials}
\setcounter{equation}{0}
The Askey-Wilson polynomials are closely related with the so-called AW(3)-algebra \cite{Zhe}, \cite{Ter}. Namely, there is a "canonical" representation of the $AW(3)$ algebra  such that one of its generators has the Askey-Wilson polynomials as eigenfunctions.

$AW(3)$ algebra usually consists of 3 generators \cite{Zhe}. Among different equivalent forms of the $AW(3)$ algebra we choose the following one, which possesses an obvious symmetry with respect to all 3 generators
(see, e.g. \cite{IT}):
\be
XY-qYX=\mu_3 Z+\omega_3, \quad YZ-qZY=\mu_1 X+\omega_1, \quad ZX-qXZ=\mu_2 Y+\omega_2 \lab{AW3} \ee
Here $q$ is a fixed parameter corresponding to the "base" parameter in $q$-hypergeometric functions for the Askey-Wilson polynomials \cite{KS}.
The pairs of operators $(X,Y)$, $(Y,Z)$ and $(Z,X)$ play the role of "Leonard pairs" (see \cite{Ter}, \cite{IT}).

The Casimir operator
\be
Q= (q^2-1)XYZ + \mu_1 X^2 + \mu_2
q^2 Y^2 + \mu_3 Z^2 +(q+1)(\omega_1 X + \omega_2 q Y + \omega_3 Z) \lab{Cas_Q} \ee
commutes with all operators $X,Y,Z$.

The constants $\omega_i, \: i=1,2,2$ (together with the value of the Casimir operator $Q$) define representations of the $AW(3)$ algebra
(see \cite{Zhe} for details).

In the case of the little $q$-Jacobi operator, the realization of the AW(3) algebra is given by the operators
\be
X= g(L+1+qab), \quad Y=x \lab{XY_q_Jac} \ee
where
$$
g=\frac{1}{(q^2-1)\sqrt{ab}}
$$
(the arithmetic meaning of the square root is assumed),
and the operator $L$ coincides with the difference eigenvalue operator for the little $q$-Jacobi polynomials in lhs of \re{lqJ_eq}, i.e.
\be
L f(x) = a(bq-x^{-1})(f(qx)-f(x))
+(1-x^{-1})(f(q^{-1}x)-f(x)) \lab{L-qJ} \ee
The operator $Z$ is defined as
$$
Zf(x) = \frac{(1-x) f(x/q)}{q \sqrt{ab}}.
$$
We then have the $AW(3)$ relations
\be
XY-qYX= Z + \omega_3, \quad YZ-qZY =0, \quad ZX-qXZ=Y+\omega_2, \lab{AW-qJ} \ee
where
$$
\omega_2=-\frac{b+1}{b(q+1)}, \quad \omega_3 = -\frac{a+1}{\sqrt{ab}(q+1)}.
$$
The Casimir operator
\be
Q=(q^2-1)XYZ  +
q^2 Y^2 + Z^2 +(q+1)(\omega_2 q Y + \omega_3 Z) \lab{Cas_lqJ} \ee
takes the value
$$
Q=-b^{-1}.
$$
These relations survive in the limit $q=-1$.

Indeed, let us define the operators
\be
X= \frac{L_0}{2} -\frac{1+\alpha+\beta}{2}, \quad Y=x, \quad Z=(x-1)R,
\lab{X_op} \ee
where $L_0$ is the operator defined by \re{L_0_R} and the operator $Y$ is multiplication by $x$.

Then it is elementary to verify that the operators $X,Y,Z$
satisfy the relations \be XY+YX=Z+\alpha, \quad YZ+ZY=0, \quad  ZX+XZ=Y+\beta
 \lab{XYZ_alg} \ee
which corresponds to the $AW(3)$ algebra with parameters $q=-1, \: \omega_3=\alpha, \: \omega_1=0, \: \omega_2=\beta$.

It is easily verified that the Casimir operator commuting with
$X,Y,Z$ is \be Q=Y^2+Z^2. \lab{Casimir_Q} \ee  In the case of of the realization \re{X_op} of the operators $X,Y,Z$,
the Casimir operator becomes the identity operator: \be Q= I.
\lab{Q_value} \ee
Note that relations \re{XYZ_alg} can be considered as an anticommutator version of some Lie algebra. Like in the case of an ordinary Lie algebra,
the Casimir operator is quadratic in the operators $Y,Z$ (the cubic part $XYZ$ disappears in both ``classical`` limits $q=\pm 1$).

The ''canonical`` representation of the algebra \re{XYZ_alg} is obtained in the basis $P_n(x)$ of orthogonal polynomials.
In this basis the operator $X$ is diagonal
$$
X P_n(x) = \frac{\lambda_n + 1+\alpha+\beta}{2} P_n(x)
$$
while the operator $Y$ is 3-diagonal
$$
Y P_n(x) = P_{n+1}(x) + b_n P_n(x) + u_n P_{n-1}(x).
$$
The monomial basis $M_n=x^n$ provides another convenient representation, where the lower- and upper- triangular operators $X,Y$ are
$$
X M_n = \xi_n M_n + \eta_n M_{n-1}, \quad Y M_n = M_{n+1},
$$
with the coefficients $\xi_n, \: \eta_n$ given by \re{xi_eta_lim} .

\section{The special case $\alpha=0$ and a "square root" of the Schr\"odinger operator}
\setcounter{equation}{0}
Consider the special case $\alpha=0$ of the little -1 Jacobi polynomials. From the expression for the weight function \re{weight} it is seen that
\be
w(x) = \kappa (1-x^2)^{(\beta-1)/2}(1+x).
\lab{w_clas} \ee
The weight function \re{w_clas} corresponds to the weight of the classical Jacobi polynomials \cite{KS} $P_n^{(a,a+1)}(x)$, where $a=(\beta-1)/2$.
Moreover, from the Dunkl-classical property \re{TPP} it is seen that for $\alpha=0$ the Dunkl operator $T_{\mu}$ becomes the ordinary derivative operator and we obtain the following realization of the Hahn property
\be
\partial_x P_n^{(a,a+1)}(x)= n P_{n-1}^{(a+1,a+2)}(x),
\lab{Hahn_Jac} \ee
which characterizes classical orthogonal polynomials \cite{Ismail_book}. Thus, in the special case $\alpha=0$, the little -1 Jacobi polynomials become the ordinary Jacobi polynomials $P_n^{(a,a+1)}(x)$. For the sake of brevity, in what follows we will denote $Q_n(x) = P_n^{(a,a+1)}(x)$, i.e. $Q_n(x)$ are the little -1 Jacobi polynomials with $\alpha=0$ and $\beta=2a+1$.

It is well known, on the one hand, that the ordinary Jacobi $P_n^{(a,b)}(x)$ polynomials satisfy the Sturm-Liouville eigenvalue equation \cite{KS}
\be
S P_n^{(a,b)}(x) = -\lambda_n \: P_n^{(a,b)}(x),
\lab{Sturm_Jacobi} \ee
where $\lambda_n=n(n+a+b+1)$ and
\be
S= (1-x^2) \partial_x^2 +\left(b-a -(a+b+2)x \right) \partial_x. \lab{L_ordinary} \ee
In the special case $b=a+1$, we have the operator
\be
S_0= (1-x^2) \partial_x^2 +\left(1 -(2a+3)x \right)\partial_x.  \lab{L_aa} \ee
On the other hand, we also know that the polynomials $Q_n(x)$ (being a special case of the little -1 Jacobi polynomials) satisfy the eigenvalue equation \re{L_0_P}
\be
L_0 Q_n(x) = \lambda_n Q_n(x), \lab{spec_eig} \ee
where the operator $L_0$ now is
\be
L_0= 2(1-x) \partial_x  R + 2(a+1) (1-R). \lab{L_0_spec} \ee
We thus obtain a rather surprising result: the ordinary Jacobi polynomials $P_n^{(a,a+1)}(x)$ satisfy not only the Sturm-Liouville equation \re{Sturm_Jacobi} but also the Dunkl-type eigenvalue equation \re{spec_eig}.

Further analysis shows however that the operator $L_0$ commutes with the operator $S_0$ defined by \re{L_aa}:
$$
L_0 S_0 - S_0 L_0=0
$$
and that there is a simple algebraic relation between these operators:
\be
L_0^2 - 4(1+a) L_0 =-4S_0. \lab{L0S0} \ee
Burchnall and Chaundy developed the theory of commuting ordinary differential operators $V_1$ and $V_2$
(the order of the operators may be arbitrary) \cite{BC}. They showed that if these operators commute, i.e.  $V_1V_2=V_2 V_1$ then necessarily they satisfy an algebraic equation (see also \cite{Ince})
\be
{\cal P}(V_1,V_2)=0, \lab{Q=0} \ee
where ${\cal P}(V_1,V_2)$ is a polynomial of two such commuting operators.

In our case we have a nontrivial example of the two commuting operators: $S_0$ is the ordinary second-order Sturm-Liouville differential operator while $L_0$ is a first-order differential operator of Dunkl type. These operator satisfy an algebraic equation \re{L0S0} as well.

There is also an interesting quantum mechanical interpretation of this example.

Let us perform the change of independent variable  $x=\sin y$ in the operator $S_0$:
$$
S_0 = \partial_y^2 + \left( \frac{1-2(a+1)\sin y}{\cos^2 y} \right) \partial_y.
$$
Using the function
\be
\Phi(y) = \sqrt{1+\sin y} \: \cos^{a+1/2} y,
\lab{Phi_y} \ee
we can transform this operator to the Schr\"odinger operator
\be
H \equiv  \Phi(y)S_0 \Phi^{-1}(y) = \partial_y^2 - U(y) \lab{Schr_H} \ee
with the potential
\be
U(y) = -(a+1)^2 + \frac{(a+1/2)(a+1/2 - \sin y)}{\cos^2 y}. \lab{U_y} \ee
This potential belongs to the class of generalized P\"oschl-Teller potentials and the corresponding Schr\"odinger equation can be solved explicitly (see, e.g. \cite{Cooper}).

Under the change of variable $x=\sin y$ and the similarity transformation $L_0 \to \Phi(y) L_0 \Phi^{-1}(y)$ the operator $L_0$ becomes
\be
L_0 =  2 \left(\partial_y -\frac{a+1/2}{\cos y} \right) R_y + 2(a+1),  \lab{L0y} \ee
where $R_y f(y)=f(-y)$ is the reflection operator with respect to the variable $y$. The constant term $2(a+1)$ is inessential and we can introduce operator
\be
 L_1 =   \left(\partial_y -\frac{a+1/2}{\cos y} \right) R_y
 \lab{L1y} \ee
without this term. Similarly, we can define the Hamiltonian
\be
H_1 = -\partial_y^2 + \frac{(a+1/2)(a+1/2 - \sin y)}{\cos^2 y} \lab{H_1} \ee
without  the superfluous constant term. (We changed the sign in the front of \re{Schr_H} in order to adjust to the canonical definition of a Schr\"odinger Hamiltonian \cite{Cooper}).

It is then easily  verified that the operator $L_1$ commutes with the Hamiltonian $H_1$: $L_1 H_1 = H_1 L_1$ and that
\be
L_1^2 =  H_1. \lab{L1H1} \ee
Relation \re{L1H1} shows that the operator $L_1$ defined by \re{L1y} can be considered as a "square root" of the Schr\"odinger Hamiltonian \re{H_1}.

In what follows we will assume that $a>1/2$. The potential  $U(y)$ of the Schr\"odinger Hamiltonian  \re{H_1} is an infinitely deep well with walls at $x = \pm \pi/2$. Clearly, $\lim_{y \to- \pm \pi/2} U(y) = \+ \infty$. Hence the spectrum of the corresponding Hamiltonian is purely discrete and infinite.
The corresponding wave functions $\psi_n(y)$ are square integrable on the interval $[-\pi/2,\pi/2]$ and are solutions of the Schr\"odinger equation
$$
H_1 \psi_n(y) =E_n \psi_n(y), \quad n=0,1,2,\dots
$$
with the energy
\be
E_n = (a+n+1)^2. \lab{E_n} \ee
Explicitly the wave functions are
\be
\psi_n(y)= \psi_0(y) \: {_2}F_1 \left( {-n, n+2a+2 \atop a+1} \left |  \frac{1-\sin y}{2} \right .  \right), \lab{psi_n} \ee
where
\be
\psi_0(y) = \Phi(y) = \sqrt{1+\sin y} \: \cos^{a+1/2} y \lab{vacuum_psi} \ee
is wave function of the ground state. Note that $\psi(\pm \pi/2)=0$ as required for eigenvalue solutions in the infinitely deep well.

The operator $L_1$ acts on the eigenfunctions as
\be
L_1 \psi_n(y) = (-1)^{n+1} (a+n+1) \psi_n(y) \lab{L_1_psi} \ee
which corresponds to relation \re{L1H1}.

This example can be generalized in the following way.

Consider the Schr\"odinger Hamiltonian
\be
H= -\partial_y^2 +U(y) \lab{Sch_H} \ee
with a potential $U(y)$, and the Dunkl-type operator
\be
L= (\partial+\chi(y))R \lab{L_Sch} \ee
with some function $\chi(y)$.

When do the operators $H$ and $L$ commute? It is easily verified that the commutativity condition $HL=LH$ is equivalent to the conditions
\be
2 \chi'(y) = U(y)-U(-y) \lab{chi_U1} \ee
and
\be
\chi''(y)  = u'(-y) + \chi(y) (U(y)-U(-y)). \lab{chi_U2} \ee
From condition \re{chi_U1} it follows that $\chi'(y)$ is an odd function. If we assume that $\chi(y)$ has a finite derivative for $y=0$ then necessarily $\chi(y)$ should be an even function, i.e. $\chi(-y)=\chi(y)$.

Condition \re{chi_U2} can then be rewritten as
\be
\chi''(y)  = U'(-y) + \chi(y) (U(y)-U(-y)) = U'(-y) + 2 \chi(y) \chi'(y) \lab{chi_U3} \ee
from which we find
\be
\chi'(y) = -U(-y) + \chi^2(y) - C \lab{cond_chi_U1} \ee
with some constant $C$. Changing sign in $y$ we have similarly
\be
-\chi'(y) = -U(y) + \chi^2(y) - C \lab{cond_chi_U2} \ee
Subtracting \re{cond_chi_U1} and \re{cond_chi_U2} we arrive at \re{chi_U1} which means that conditions \re{cond_chi_U1} and \re{cond_chi_U2} are compatible with \re{chi_U1} and \re{chi_U2}. Adding \re{cond_chi_U1} and \re{cond_chi_U2}, we have
\be
2 \chi^2 (y) = U(y)+U(-y) + 2C. \lab{cond_chi_U3} \ee
There is a remarkable connection with the factorization method of the Schr\"odinger equation \cite{Cooper}, \cite{Spi}. Indeed, assume that for the potential $V(y)$ and a constant $C$ there is a the representation
\be
H +C \equiv  -\partial_y^2 + V(y) +C = A^{\dagger} A \lab{factr_H} \ee
where
\be
A=\partial_y + \chi(y), \quad A^{\dagger}=-\partial_y + \chi(y)
\lab{AA_def} \ee
The Schr\"odinger Hamiltonian H is then said to be factorizable in terms of the operators $A^{\dagger}, A$.

The potential $V(y)$ is related with with the "superpotential" $\chi(y)$ by obvious fromula
\be
\chi^2(y) -\chi'(y) = V(y) +C \lab{V_chi} \ee
A refactorization performed by permuting the operators $A^{\dagger}, A$ leads to the new Hamiltonian
\be
\t H +C \equiv -\partial_y^2 + \t V(y) + C = A A^{\dagger},
\lab{tHf} \ee
where the new potential $\t V(y)$ is defined as
\be
\chi^2(y) + \chi'(y) = \t V(y) +C. \lab{tV_def} \ee
Comparing formulas \re{V_chi} and \re{tV_def} with \re{cond_chi_U1} and \re{cond_chi_U2}, we see that one can put
$$
V(y) = U(-y), \quad \t V(y) = U(y),
$$
i.e. the refactorization procedure simply leads to a change of sign in the argument of the potential: $U(-y) \to U(y)$. Such potentials were studied by Spiridonov and are special cases of more general self-similar potentials (see \cite{Spi}, \cite{Spi_PR} for details).

We thus see that for any smooth even function $\chi(y)$ the Dunkl-type operator $L$ \re{L_Sch} commutes with the Hamiltonian $H$ \re{Sch_H}, where the potential $U(y)$ is defined by \re{cond_chi_U2} (the constant $C$ can be arbitrary, of course).

It is simply verified that in this case we have
\be
L^2 = (\partial_y + \chi(y))R(\partial_y + \chi(y))R = (\partial_y + \chi(y))(-\partial_y + \chi(y))=-\partial_y^2 +U(y) +C = H+C \lab{L_2_H} \ee
i.e. the commuting operators $L$ and $H$ satisfy a simple algebraic relation of the second order.

Thus one can say that the operator $L$ is a "square root" of the Schr\"odinger operator for potentials $U(y)$. Of course, such operator exists only for potentials $U(y)$ satisfying condition \re{cond_chi_U2} with some even function $\chi(y)$.

The physical meaning of the operator $L$ is very simple. Indeed, let $\psi(x)$ be an eigenfunction of the Hamiltonian \re{Sch_H}
\be H \psi(y)  = E \psi(y) \lab{HpsiE} \ee
corresponding to the energy $E$. Obviously,
\be
L \psi(y) = (\pm) \sqrt{E} \psi(y). \lab{LpsiE} \ee
Acting on \re{LpsiE} with the operator $R$ we get
\be
RL \psi(y)  = A^{\dagger} \psi(y) = (\pm) \sqrt{E} \psi(-y), \lab{RLpsi} \ee
where the operator $A^{\dagger}$ is defined by \re{AA_def}. We see that the operator $RL$ is the Darboux operator which transforms the eigenfunction $\psi(y)$ into the function $\psi(-y)$ which is the eigenfunction of the Hamiltonian
$$
\t H = -\partial_y^2 + U(-y).$$
Thus the operator $RL$ "flips" the eigenfunctions corresponding to the potentials $U(y)$ and $U(-y)$.

\bigskip\bigskip\bigskip
{\Large\bf Acknowledgments}
\bigskip

\noindent The authors are indebted to R.Askey, C.Dunkl, T.
Koornwinder, W.Miller, V.Spiridonov, P.Terwilliger and P.Winternitz for stimulating communications
and to a referee for many useful remarks and for drawing their attention to formula \re{f_elem}. AZ
thanks CRM (U de Montr\'eal) for its hospitality.

\bigskip\bigskip

\bb{99}


\bi{AI} R. A. Askey and M. E. H. Ismail, {\it Recurrence relations, continued fractions and orthogonal polynomials}, Mem. AMS, {\bf 49} (1984), No. 300,  1--108.

\bi{Atia} M. J. Atia, {\it An example of non-symmetric
semi-classical form of class s=1. Generalization of a case of
Jacobi sequence}. Int. J. Math. and Math. Sci. {\bf 24} (10)
(2000), 673--689.

\bi{Atia2} M.J.Atia, {\it Some~generalized~Jacobi~polynomials},
OPSFA 2009, Proceedings (to be published). 


\bi{BI} E. Bannai and T. Ito, {\it Algebraic Combinatorics I: Association Schemes}. 1984. Benjamin \& Cummings, Mento Park, CA.

\bi{Belm} S. Belmehdi, {\it Generalized Gegenbauer polynomials}, J. Comput. Appl Math. {\bf 133} (2001), 195-–205.

\bi{Cheikh} Y. Ben Cheikh and M.Gaied, {\it Characterization of the Dunkl-classical symmetric orthogonal polynomials}, Appl. Math. and Comput.
{\bf 187}, (2007) 105--114.

\bi{BC} J. L. Burchnall,T. W. Chaundy, {\it Commutative ordinary differential operators}, Proc.
London Math. Soc. (Ser. 2), {\bf 21} (1922),420�-440.


\bi{Chi} T. Chihara, {\it An Introduction to Orthogonal
Polynomials}, Gordon and Breach, NY, 1978.

\bi{Chi_Chi} L.M.Chihara,T.S.Chihara, {\it A class of nonsymmetric orthogonal polynomials}. J. Math. Anal. Appl. {\bf 126} (1987), 275--291.

\bi{Cooper} Cooper F., Khare A., Sukhatme U., {\it Supersymmetry and quantum mechanics}, Phys. Rep. {\bf 251} (1995), 267-�
385, hep-th/9405029.

\bi{Dunkl} C.F.Dunkl, {\it Integral kernels with reflection group invariance}. Canadian Journal of Mathematics, {\bf 43} (1991)
1213--1227.

\bi{Dunkl2} C.F.Dunkl, {\it Orthogonal polynomials and reflection groups}, in: Special functions 2000: current perspective and future directions, NATO Science Series, 111--128, 2001.


\bi{Ger1} Ya.L.Geronimus, {\it On polynomials orthogonal with respect to to
the given numerical sequence and on Hahn's theorem}, Izv.Akad.Nauk, {\bf 4}
(1940), 215-228 (in Russian).







\bi{Ince} E.L. Ince, {\it Ordinary Differential Equations}, New York, Dover, 1956.

\bi{Ismail_book} M.E.H.Ismail, {\it Classical and Quantum orthogonal polynomials in one variable}.
Encyclopedia of Mathematics and its Applications (No. 98), Cambridge, 2005.

\bi{IT} T.Ito, P.Terwilliger, {\it Double Affine Hecke Algebras of Rank 1
and the $Z_3$-Symmetric Askey-Wilson Relations}, SIGMA {\bf 6} (2010), 065, 9 pages. arXiv:1001.2764.

\bi{KS} R.Koekoek, R.Swarttouw, {\it The Askey-scheme of hypergeometric orthogonal polynomials and its $q$-analogue}, Report no. 98-17, Delft University of Technology, 1998.

\bi{KLS} R. Koekoek,P. Lesky, R. Swarttouw, {\it Hypergeometric Orthogonal Polynomials and Their $Q$-analogues}, Springer-Verlag, 2010.



\bi{Leonard} D.Leonard, {\it Orthogonal Polynomials, Duality and Association Schemes}, SIAM J. Math. Anal. {\bf 13} (1982) 656--663.






\bi{Spi_PR} V.Spiridonov, {\it Universal Superpositions of Coherent States
and Self-Similar Potentials}, Phys. Rev. {\bf A52}, (1995), 1909--1935.

\bi{Spi} V.Spiridonov, {\it The factorization method, self-similar
potentials and quantum algebras}, in: �Special Functions 2000:
Current Perspective and Future Directions�, Proc. NATO Advanced
Study Institute (Tempe, USA, May 29 � June 9, 2000), edited by
J. Bustoz, M.E.H. Ismail and S.K. Suslov, Kluwer Academic
Publishers, Dordrecht, 2001, pp. 335-�364,
arXiv:hep-th/0302046v1.

\bi{SVZ} V.Spiridonov, L.Vinet and A.Zhedanov, {\it Spectral
transformations, self-similar reductions and orthogonal polynomials}, J.Phys.
A:  Math.  and Gen.  {\bf 30} (1997), 7621--7637.

\bi{Sz} G. Szeg\H{o}, Orthogonal Polynomials, fourth edition,  AMS, 1975.

\bi{Ter} P. Terwilliger, {\it Two linear transformations each tridiagonal with respect to an eigenbasis of the other}.
Linear Algebra Appl. {\bf 330} (2001) 149--203.




\bi{Zhe} A. S. Zhedanov. {\it ''Hidden symmetry`` of Askey-Wilson polynomials}, Teoret. Mat. Fiz.
{\bf 89} (1991) 190--204. (English transl.: Theoret. and Math. Phys. {\bf 89} (1991), 1146--1157).


\bibitem{ZhS} A.S. Zhedanov, {\it Rational spectral transformations
and orthogonal polynomials}, J. Comput. Appl. Math. {\bf 85}, no. 1
(1997), 67--86.

\eb

\end{document}